\ifx\documentclass\undefined
\documentstyle[12pt]{article}
\else
\documentclass[12pt]{article}
\usepackage{latexsym}
\usepackage{amsfonts}
\usepackage{amsmath}
\usepackage{amssymb}
\fi

\author{ M\'at\'e Bezdek}

\font\tenBbb=msbm10 at 12pt         \font\sevenBbb=msbm9    \font\fiveBbb=msbm7
\newfam\Bbbfam
\textfont\Bbbfam=\tenBbb \scriptfont\Bbbfam=\sevenBbb
\scriptscriptfont\Bbbfam=\fiveBbb

\def\E{{\mathbb E}}

\def\kkk{\null\hfill $\Box$\smallskip}

\def\HH{{\bf H}}

\def\CC{{\bf C}}
\def\DD{{\bf D}}

\newcommand{\proof}{{\noindent\bf Proof:{\ \ }}}

\newtheorem{theorem}{Theorem}[section]
\newtheorem{lemma}[theorem]{Lemma}
\newtheorem{sublemma}[theorem]{Sublemma}
\newtheorem{remark}[theorem]{Remark}

\newtheorem{cor}[theorem]{Corollary}

\newtheorem{definition}[theorem]{Definition}

\newtheorem{prob}[theorem]{Problem}

\title{On a generalization of the Blaschke-Lebesgue theorem for disk-polygons
}

\begin{document}

\maketitle

\date

\begin{abstract}
In this paper we prove an extension of the Blaschke-Lebesgue theorem for a family of convex domains called disk-polygons. Also, this provides yet another new proof of the Blaschke-Lebesgue theorem.
\end{abstract}

\section{Introduction}
\label{zero}

A {\it convex domain} of the Euclidean plane $\E^2$ is a compact convex set with non-empty interior. Let $\CC\subset\E^2$ be a convex domain, and let $l\subset\E^2$ be a line. Then the distance between the two supporting lines of $\CC$ parallel to $l$ is called the width of $\CC$, in {\it direction} $l$. Moreover, the smallest width of $\CC$ is called the {\it minimal width} of $\CC$, labelled by $w(\CC)$. In other words, the minimal width of a convex domain is equal to the smallest distance between parallel supporting lines of the given convex domain. Also, recall that the convex domain $\CC\subset\E^2$ is called a {\it convex domain of constant width} $w$, if the width of $\CC$ in any direction of $\E^2$ is equal to $w$. The simplest example of a convex domain of constant width $w$ is the circular disk of diameter $w$. However, the family of convex domains of constant width $w$ is a large and rather complex family. For example, a {\it Reuleaux polygon} of width $w$ is a convex domain of constant width $w$, whose boundary is a union of finitely many circular arcs of radii $w$. (For a detailed account on a number of elementary properties of convex domains of constant width see for example \cite{BoYa61}.) The simpliest example of a Reuleaux polygon is the {\it Reuleaux triangle}. To construct a Reuleaux triangle of width $w$, start with an equilateral triangle of side length $w$; then take the intersection of the three circular disks of radii $w$, centered at the vertices of the equilateral triangle. In fact, the family of Reuleaux polygons of width $w$ is a dense subset of the family of convex domains of constant width $w$. (For more details on this see \cite{BoYa61}.) Perhaps, it is then not surprising that there are convex domains of constant width whose boundaries include no circular arcs, however small. (For a very flexible way of constructing convex domains of constant width see \cite{Sa70}.)

On the one hand, the classical isoperimetric inequality combined with Barbier's theorem (stating that the perimeter of any convex domain of constant width $w$ is equal to $\pi w$) implies that the largest area of convex domains of constant width $w$, is the circular disk of diameter $w$, having the area of $\frac{\pi}{4}w^2$ (for more details see for example \cite{BoYa61}). On the other hand, the well-known Blaschke-Lebesgue theorem states that among all convex domains of constant width $w$, the Reuleaux triangle of width $w$ has the smallest area, namely $\frac{1}{2}(\pi-\sqrt{3})w^2$. W. Blaschke \cite{Bla15} and H. Lebesgue \cite{Leb14} were the first to show this and the succeeding decades have seen other works published on different proofs of that theorem. For a most recent new proof, and for a survey on the state of the art of different proofs of the Blaschke-Lebesgue theorem, see the elegant paper of E. M. Harrell \cite{Ha02}.

The main goal of this paper is to provide yet another new proof of the Blaschke-Lebesgue theorem, and perhaps, more importantly, to prove a new more general version of it with the hope of extending it to higher dimensions. In the remaining part of the introduction we summarize our new results starting with the necessary definitions.

Our first definition has been introduced in \cite{Be07} and it specifies the type of sets studied in this paper.

\begin{definition}
{\rm The intersection of finitely many (closed) circular disks of unit radii with non-empty interior in $\E^2$ is called a {\it disk-polygon}. We will assume that whenever we take a disk-polygon, then the disks generating it, simply called {\it generating disks}, are all needed; that is, each of them contributes to the boundary of the disk-polygon through a circular arc called a {\it side}, with the consecutive pairs of sides meeting in the {\it vertices} of the given disk-polygon.} 
\end{definition}

The parameter introduced in the next definition turns out to be a crucial one for our investigations.

\begin{definition}
{\rm The disk-polygon $\DD$ is called a disk-polygon with {\it center parameter} $d$, $0 < d < \sqrt{3} = 1.732\dots$, if the distance between any two centers of the generating disks of $\DD$ is at most $d$. Let ${\cal{F}}(d)$ denote the family of all disk-polygons with center parameter $d$. }
\end{definition}

The following special disk-polygon is going to play a central role in our investigations. 

\begin{definition}
{\rm Let  $\Delta(d)$ denote the regular disk-triangle whose three generating (unit) disks are centered at the vertices of a regular triangle of side length $d$, $1\le d < \sqrt{3} = 1.732\dots$.}
\end{definition}  

Recall, that the {\it inradius} $r(\CC)$ of a convex domain $\CC$ in $\E^2$ is the radius of the largest circular disk lying in $\CC$ (simply called the {\it incircle} of $\CC$). 

\begin{remark}
{\rm The following formulas give  the inradius $r(\Delta(d))$, the minimal width $w(\Delta(d))$, the area $a(\Delta(d))$ and the perimeter $p(\Delta(d))$ of $\Delta(d)$ for all $1\le d < \sqrt{3}$:

$$r(\Delta(d))=1-\frac{1}{3}\sqrt{3}d;$$

$$w(\Delta(d))=1-\frac{1}{2}\sqrt{4+2d^2-2\sqrt{3}d\sqrt{4-d^2}};$$

$$a(\Delta(d))=\frac{3}{2}\arccos d+\frac{1}{4}\sqrt{3}d^2-\frac{3}{4}d\sqrt{4-d^2}-\frac{1}{2}\pi;$$

$$p(\Delta(d))=2\pi-6\arcsin \frac{d}{2}.$$
}
\end{remark}

Now, we are ready to state our first theorem.

\begin{theorem}\label{elso}
Let $\DD \in {\cal{F}}(d)$ be an arbitrary disk-polygon with center parameter $d, 1 \le d < \sqrt{3}$. Then, the area of $\DD$ is at least as large as the area of $\Delta(d)$,  i.e. 
$$a(\DD)\ge a(\Delta(d))$$
with equality if and only if $\DD=\Delta(d)$.
\end{theorem}

\begin{remark}
{\rm For $d=1$ the above area inequality and the well-known fact (see for example \cite{BoYa61}) that the family of Reuleaux polygons of width $1$ is a dense subset of the family of convex domains of constant width $1$, imply the Blaschke - Lebesgue theorem in a straighforward way. }
\end{remark}

In connection with Theorem~\ref{elso} K. Bezdek \cite{Be08} proposed to investigate the following related problem.

\begin{prob}
Let $\DD \in {\cal{F}}(d)$ be an arbitrary disk-polygon with center parameter $d, 1 \le d < \sqrt{3}$. Prove or disprove that the perimeter of $\DD$ is at least as large as the perimeter of $\Delta(d)$,  i.e. 
$$p(\DD)\ge p(\Delta(d)).$$
\end{prob}

As the following two statements belong to the core part of our proof of Theorem~\ref{elso} and might be of independent interest, we mention them here.

\begin{lemma}\label{first}
Let $\DD \in {\cal{F}}(d)$ be an arbitrary disk-polygon with center parameter $d, 1 \le d < \sqrt{3}$. Then, the inradius of $\DD$ is at least as large as the inradius of $\Delta(d)$,  i.e. 
$$r(\DD)\ge r(\Delta(d) ).$$
\end{lemma}

\begin{lemma}\label{second}
Let $\DD \in {\cal{F}}(d)$ be an arbitrary disk-polygon with center parameter $d, 1 \le d < \sqrt{3}$. Then, the minimal width of $\DD$ is at least as large as the minimal width of $\Delta(d)$,  i.e. 
$$w(\DD)\ge w(\Delta(d) ).$$
\end{lemma}

Let $\CC\subset\E^2$ be a convex domain and let $\rho>0$ be given. Then, the {\it outer parallel domain} $\CC_{\rho}$ of radius $\rho$ of $\CC$ is the union of all (closed) circular disks of radii $\rho$, whose centers belong to $\CC$. Recall that $a(\CC_{\rho}) = a(\CC) + p(\CC){\rho} + \pi{\rho}^2$.

\begin{definition}
{\rm Let  $0 < d < 1$ be given and let $\Delta^{\circ}(d)$ denote the outer parallel domain of radius $1-d$ of a Reuleaux triangle of width $d$.} 
\end{definition}

\begin{remark}
{\rm Note that $\Delta^{\circ}(d)$ is a convex domain of constant width $2-d$ and so, Barbier's theorem (\cite{BoYa61}) implies that its perimeter is equal to $p(\Delta^{\circ}(d))=\pi (2-d)$ moreover, it is not hard to check that its area is equal to $a(\Delta^{\circ}(d))=\frac{1}{2}(\pi-\sqrt{3})d^2-\pi d+\pi $.}
\end{remark}

Now, we are ready to state our second theorem.

\begin{theorem}\label{masodik}
Let $\DD \in {\cal{F}}(d)$ be an arbitrary disk-polygon with center parameter $d, 0<d<1$. Then, the area of $\DD$ is strictly larger than the area of $\Delta^{\circ}(d)$,  i.e. 
$$a(\DD)> a(\Delta^{\circ}(d)).$$
\end{theorem}

\begin{remark}
{\rm Note that our proof of Theorem~\ref{masodik} implies that the above lower bound is best possible.} 
\end{remark}

In connection with Theorem~\ref{masodik} K. Bezdek \cite{Be08} has raised the following question.

\begin{prob}
Let $\DD \in {\cal{F}}(d)$ be an arbitrary disk-polygon with center parameter $d, 0<d<1$. Prove or disprove that the perimeter of $\DD$ is strictly larger than the perimeter of $\Delta^{\circ}(d)$,  i.e. 
$$p(\DD)> p(\Delta^{\circ}(d)).$$
\end{prob}

\medskip
\section{Proof of Lemma~\ref{first}}

The following definition introduces the notion of {\it dual disk-polygon} that turns out to play a central role in our investigations.

\begin{definition}
{\rm Let $\DD$ be an arbitrary disk-polygon in $\E^2$. Then the intersection of the circular disks of unit radii centered at the vertices of $\DD$ is called the {\it dual disk-polygon} $\DD^*$ associated with $\DD$.} 
\end{definition}

\begin{remark}\label{fontos}
{\rm It is easy to see that $(\DD^*)^*=\DD$ for any disk-polygon $\DD\subset \E^2$.}
\end{remark}

For the sake of completeness we recall the following definition as well. (For more details see \cite{Be07}.)

\begin{definition} 
{\rm Let ${\bf X}$ be an arbitrary set contained in a unit circular disk of $\E^2$. Then the {\it spindle convex hull} of ${\bf X}$ is the intersection of all the unit circular disks that contain ${\bf X}$. Moreover, we say that ${\bf X}$ is {\it spindle convex} if for any two points of ${\bf X}$ their spindle convex hull is contained in ${\bf X}$.}
\end{definition}

\begin{remark}
{\rm It is easy to see that if $\DD$ is an arbitrary disk-polygon in $\E^2$, then its dual disk-polygon $\DD^*$ is the spindle convex hull of the centers of the generating disks of $\DD$. Moreover, $\DD$ is spindle convex.}
\end{remark}

Recall, that the {\it circumradius} $R(\CC)$ of a convex domain $\CC$ in $\E^2$ is the radius of the smallest circular disk containing $\CC$ (simply called the {\it circumcircle} of $\CC$).

\begin{sublemma}\label{sublemma1} 
If $\DD$ is an arbitrary disk-polygon in $\E^2$, then the incircle of $\DD$ and the circumcircle of $\DD^*$ are concentric circular disks moreover,
$$r(\DD)+R(\DD^*)=1.$$
\end{sublemma}

\proof
Let the incircle of $\DD$ be centered at $C$ having radius $r(\DD)$. Note that Remark~\ref{fontos} implies that the centers of the circular disks generating the sides of $\DD$ are identical to the vertices of $\DD^*$, and vica versa the circular disk of unit radius centered at a vertex of $\DD^*$ generates a side for $\DD$. As a result the circular disk of radius $1-r(\DD)$ centered at $C$ contains all the vertices of $\DD^*$ and so, it contains $\DD^*$, proving that $R(\DD^*)\le 1- r(\DD)$. In the same way, one can show that starting with the circumcircle of $\DD^*$ centered say, at $C^*$ having radius $R(\DD^*)$, the circular disk centered at $C^*$ having radius $1-R(\DD^*)$ is contained in $\DD$ and therefore $1-R(\DD^*)\le r(\DD)$, finishing the proof of Sublemma~\ref{sublemma1}.
\kkk

Clearly, Sublemma~\ref{sublemma1} implies in a straighforward way that the inequality of Lemma~\ref{first} is equivalent to the following one:
$$R(\DD^*)\le R(\Delta(d)^*)=\frac{d}{\sqrt{3}}<1.$$

Now, recall that according to the well-known Jung theorem (see for example \cite{BoYa61}) the circumradius of a finite point set of diameter $d$ in $\E^2$ is at most $\frac{d}{\sqrt{3}}$ (which is in fact, the circumradius of a regular triangle of side length $d$). Also, note that the smallest circular disk containing $\DD^*$ is identical to the smallest circular disk containing the vertices of $\DD^*$. Thus, as $\DD \in {\cal{F}}(d)$ i.e. the parwise distances between the vertices of $\DD^*$ are at most $d$, therefore Jung's theorem implies in a straighforward way that $R(\DD^*)\le \frac{d}{\sqrt{3}}$, finishing the proof Lemma~\ref{first}.

\medskip
\section{Proof of Lemma~\ref{second}}

First, recall the following statement from \cite{BeCC06}.

\begin{sublemma}\label{sublemma2}
If $\DD$ is an arbitrary disk-polygon in $\E^2$, then the Minkowski sum of $\DD$ and its dual $\DD^*$ is a convex domain of constant width $2$. 
\end{sublemma}

Second, recall that the {\it diameter} of a set $X\subset \E^2$, denoted by ${\it diam}(X)$, is the largest distance between two points in $X$. Sublemma~\ref{sublemma2} implies the following statement in a straighforward way.

\begin{cor}\label{width-diameter}
If $\DD$ is an arbitrary disk-polygon in $\E^2$, the 
$$w(\DD)+{\it diam}(\DD^*)=w(\DD^*)+{\it diam}(\DD)=2.$$
\end{cor}

Now, we are ready to give a proof of Lemma~\ref{second}. Let $\DD \in {\cal{F}}(d)$ be an arbitrary disk-polygon with center parameter $d, 1 \le d < \sqrt{3}$. Then, our goal is to show that the minimal width of $\DD$ is at least as large as the minimal width of $\Delta(d)$,  i.e. $w(\DD)\ge w(\Delta(d) )$. By Corollary~\ref{width-diameter} it is sufficient to prove that
$${\it diam}(\DD^*)\le{\it diam}(\Delta(d)^* )=1+\frac{1}{2}\sqrt{4+2d^2-2\sqrt{3}d\sqrt{4-d^2}},$$
where $\Delta(d)^*$ is the spindle convex hull of a regular triangle of side length $d$ with $1 \le d < \sqrt{3}$.
Let the line segment $AB$ represent the diameter of $\DD^*$ and let $a$ and $b$ be the lines passing through $A$ and $B$ and being perpendicular to $AB$. Clearly, $a$ and $b$ are supporting lines of $\DD^*$. Using Sublemma~\ref{sublemma2} one can assume without loss of generality that either $A$ and $B$ are both vertices of $\DD^*$ ({\it Case (i)}) or there is a side (circular arc) $s_a$ of $\DD^*$ that is tangent to $a$ at $A$ moreover, $B$ is a vertex of $\DD^*$ ({\it Case (ii)}).

{\it Case (i)}: As $A$ and $B$ are the center points of two generating disks of $\DD$ and $\DD\in {\cal{F}}(d)$, therefore ${\it diam}(\DD^*)\le d < 1+\frac{1}{2}\sqrt{4+2d^2-2\sqrt{3}d\sqrt{4-d^2}}$ for all $1 \le d < \sqrt{3}$, finishing the proof of Lemma~\ref{second} in {\it Case (i)}.

{\it Case (ii)}: Let $C$ and $D$ denote the vertices of $\DD^*$ that are the endpoints of the side $s_a$ of $\DD^*$. Clearly, as $\DD\in {\cal{F}}(d)$ therefore the line segments $BC, CD$ and $DB$ have length at most $d$. Under this condition (for given $d$) let us maximize the length of $AB$. (Here, $a$ must be tangent to $s_a$ at $A$ and $b$ must go through $B$ such that $a$ and $b$ are parallel with the triangle BCD lying between them.) It is easy to see that the maximum must belong to the geometric configuration with $BC$ and $BD$ both having their length equal to $d$ and with $A$ being the midpoint of the circular arc $s_a$. Now, it is not yet clear what the length of $CD$ should be. At this point we know only that its length is at most $d$ or putting it equivalently, the angle $2\alpha$ at the vertex $B$ of the triangle $BCD$ is at most $\frac{\pi}{3}$. However, with a bit of computation one can actually show that in order to maximize the length of $AB$ one has to have the length of $CD$ equal to $d$ as well. The details are as follows: first a simple computation yields that the length of $AB$ is equal to $f(\alpha):= d\cos\alpha-\sqrt{1-d^2\sin^2\alpha }+1$ with $0\le \alpha\le \frac{\pi}{6}$; second using for example, MAPLE one can actually check that the maximum value of $f(\alpha)$ under the condition that $0\le \alpha\le \frac{\pi}{6}$ is always  $f(\frac{\pi}{6})=1+\frac{\sqrt{3}}{2}d-\sqrt{1-\frac{1}{4}d^2}= 1+\frac{1}{2}\sqrt{4+2d^2-2\sqrt{3}d\sqrt{4-d^2}}= {\it diam}(\Delta(d)^*)$ for all $1 \le d < \sqrt{3}$ (in fact, it turns out that $f(\alpha)$ is an increasing function of $0\le\alpha\le \frac{\pi}{6}$ for all $1 \le d < \sqrt{3}$), finishing our proof in {\it Case (ii)} and completing the proof of Lemma~\ref{second}.

\medskip
\section{Proof of Theorem~\ref{elso} }

Let $\CC(O, x)$ denote the incircle of $\DD$ centered at the point $O$ having radius $x$. Either there are two diametrically opposite points of $\CC(O, x)$ say, $A$ and $B$ at which $\CC(O, x)$ touches the boundary of $\DD$ ({\it Case I}) or there three points on the boundary of $\CC(O, x)$ say, $A, B$ and $C$ at which $\CC(O, x)$ touches the boundary of $\DD$ such that $O$ belongs to the interior of the triangle $ABC$ ({\it Case II}). In both cases we show that the area of $\DD$ is at least as large as the area of $\Delta(d)$. The details are as follows.

{\it Case I}: Let $a$ (resp., $b$) be the line passing through $A$ (resp., $B$) that is perpendicular to the line segment $AB$. Clearly, $a$ and $b$ are parallel supporting lines of $\DD$ and therefore Lemma~\ref{second} implies that the distance $2x$ between them is at least $w(\Delta(d))$. Now, on the one hand, the area of $\CC(O, x)$ is at least as large as the area of a circular disk of diameter $2x$ that is it is at least as large as $f_{\rm area}(d):=\frac{\pi}{4}w^2(\Delta(d))=\frac{\pi}{4}  \bigg[1-\frac{1}{2}\sqrt{4+2d^2-2\sqrt{3}d\sqrt{4-d^2}}\bigg]^2$ on the other hand, recall that the area of $\Delta(d)$ is equal to $g_{\rm area}(d):=a(\Delta(d))=\frac{3}{2}\arccos d+\frac{1}{4}\sqrt{3}d^2-\frac{3}{4}d\sqrt{4-d^2}-\frac{1}{2}\pi $. By graphing $f_{\rm area}(d)$ and $g_{\rm area}(d)$ as functions of $d$ it is convenient to check with MAPLE that $f_{\rm area}(d)>g_{\rm area}(d)$ for all $1 \le d < \sqrt{3}$, finishing the proof of the area inequality in Theorem~\ref{elso} for {\it Case I}.

{\it Case II}:
Lemma~\ref{first} implies that $x=r(\DD)\ge r(\Delta(d))=1-\frac{1}{3}\sqrt{3}d$. Hence, using the area estimate as in {\it Case I}, we can assume that 
$$(1) \hskip1.0cm 1-\frac{1}{3}\sqrt{3}d\le x\le  \frac{1}{2}w(\Delta(d))=\frac{1}{2}-\frac{1}{4}\sqrt{4+2d^2-2\sqrt{3}d\sqrt{4-d^2}}.$$
Let $a, b$ and $c$ be the uniquely determined supporting lines of $\DD$ passing through the points $A, B$ and $C$. Moreover, let $a'$ be the line parallel to $a$ at distance $w(\Delta(d))$ from $a$ and lying on the same side of $a$ as $\DD$. Let the lines $b'$ and $c'$ be defined in a similar way using the lines $b$ and $c$. Clearly, Lemma~\ref{second} implies that there are points $A', B'$ and $C'$ with $A'\in \DD\cap a', B'\in\DD\cap b'$ and $C'\in \DD\cap c'$. As a next step take the spindle convex hull of $A'$ and $\CC(O, x)$ and subtract from it the incircle $\CC(O, x)$ and denote the set obtained by $\HH_{A'}$ (which in fact, will look like a "cap" attached to $\CC(O, x)$). In the same way, we construct the sets $\HH_{B'}$ and $\HH_{C'}$. Clearly, the sets $\CC(O, x)$, $\HH_{A'}, \HH_{B'}$ and $\HH_{C'}$ are pairwise non-overlapping moreover, as $\DD$ is spindle convex, therefore they all lie in $\DD$ and thus,
$$(2) \hskip1.0cm a(\DD)\ge a(\CC(O, x))+a(\HH_{A'})+a(\HH_{B'})+a(\HH_{C'}).$$

Finally, let $D'$ be a point exactly at distance $w(\Delta(d))-x$ from $O$. Then take the spindle convex hull of $D'$ and $\CC(O, x)$ and remove from it the incircle $\CC(O, x)$ and denote the set obtained by $\HH_{D'}$. As the length of the line segments $A'O, B'O$ and $C'O$ are all at least $w(\Delta(d))-x$, therefore introducing the notation
$F(d, x):=a(\CC(O, x))+3a(\HH_{D'})$ we get that
 
$$(3) \hskip1.0cm a(\CC(O, x))+a(\HH_{A'})+a(\HH_{B'})+a(\HH_{C'})\ge F(d, x).$$

Standard geometric calculations yield the following formula for $F(d, x)$:

$$F(d, x)=\pi x^2+3\arccos \bigg[\frac{1+2(1-x)y(d)-y^2(d) }{2(1-x)} \bigg]$$
$$-3x^2\arccos\bigg[\frac{1-(1-x)^2-(1-x-y(d))^2 }{2(1-x)(1-x-y(d)) }\bigg]$$
$$-\frac{3}{2}\sqrt{(3-2x-y(d))(1-2x-y(d))(1-y^2(d)) },$$
where $y(d)=1-w(\Delta(d))=\frac{1}{2}\sqrt{4+2d^2-2\sqrt{3}d\sqrt{4-d^2}}$.
By applying the necessary tools of MAPLE, it turns out that for each $1 \le d < \sqrt{3}$ the function $F(d, x)$
is an increasing function of $x$ over the interval defined by $(1)$. Thus,
$$(4) \hskip1.0cm F(d, x)\ge F(d, 1-\frac{1}{3}\sqrt{3}d)=a(\Delta(d)).$$
Hence, $(2), (3)$ and $(4)$ finish the proof in {\it Case II}, and so the proof of Theorem~\ref{elso} is now complete.

\medskip
\section{Proof of Theorem~\ref{masodik}}

Let $\DD \in {\cal{F}}(d)$ be an arbitrary disk-polygon with center parameter $d, 0<d<1$. Let the diameter of $\DD^*$ be denoted by $d^*$. Applying the method of {\it Case (ii)} of the proof of Lemma~\ref{second}, it turns out, that $d^*\le d$. Thus, according to a well-known theorem (see for example \cite{BoYa61}) there exists a convex domain $\bar{\DD}$ of constant width $d$ with $\DD^*\subset\bar{\DD}$. As a result we get that
$$(5) \hskip1.0cm  \DD=(\DD^*)^*\supset \bar{\DD}_{1-d},$$
where $\bar{\DD}_{1-d}$ stands for the outer parallel domain of radius $1-d$ of $\bar{\DD}$. Note that as $\bar{\DD}$ is of constant width $d$ therefore $\bar{\DD}_{1-d}$ is of constant width $2-d$. Also, it is clear that in $(5)$ equality cannot happen (simply for the reason that $\bar{\DD}_{1-d}$ is not a disk-polygon). Thus, $(5)$ implies in a straighforward way that
$$(6) \hskip1.0cm a(\DD)> a(\bar{\DD}_{1-d})=a(\bar{\DD})+p(\bar{\DD})(1-d)+\pi (1-d)^2$$
On the one hand, Barbier's theorem (\cite{BoYa61}) implies that $p(\bar{\DD})=\pi d$ and on the other hand, the Blaschke-Lebesgue theorem or better yet Theorem~\ref{elso} yields that $a(\bar{\DD})\ge \frac{1}{2}(\pi-\sqrt{3})d^2$. Thus, these facts together with $(6)$ imply that $a(\DD)>\frac{1}{2}(\pi-\sqrt{3})d^2-\pi d+\pi =
a(\Delta^{\circ}(d))$, finishing the proof of Theorem~\ref{masodik}.

\bigskip

{\bf Acknowledgement.} The author is grateful to his father K\'aroly Bezdek, for bringing the problems of this paper to his attention; and also, for the extensive discussions on the topics mentioned.

\vspace{1cm}

\medskip

\noindent
M\'at\'e Bezdek,
Notre Dame High School,
11900 Country Village Link NE,
Calgary, AB, Canada, T3K 6E4.
\newline
{\sf e-mail: mate.bezdek@hotmail.com}

\end{document}